\documentclass[11pt]{amsart2000}
\numberwithin{equation}{section}
\theoremstyle{plain}
\newtheorem{thm}{Theorem}

\theoremstyle{definition}
\newtheorem{defn}{Definition}
\newtheorem{exmp}{Example}
\newtheorem{ques}{Question}
\newtheorem*{dmc}{Double Midset Conjecture}
\newtheorem*{smc}{$(n-1)$-sphere Midset Conjecture}
\theoremstyle{remark}

\def\Bd{\text{Bd}\:}
\def\Max{\text{Max}\:}
\def\ump{\operatorname{ump}}

\begin{document}
\setlength{\unitlength}{0.01in}
\linethickness{0.01in}
\begin{center}
\begin{picture}(474,66)(0,0)
\multiput(0,66)(1,0){40}{\line(0,-1){24}}
\multiput(43,65)(1,-1){24}{\line(0,-1){40}}
\multiput(1,39)(1,-1){40}{\line(1,0){24}}
\multiput(70,2)(1,1){24}{\line(0,1){40}}
\multiput(72,0)(1,1){24}{\line(1,0){40}}
\multiput(97,66)(1,0){40}{\line(0,-1){40}}
\put(143,66){\makebox(0,0)[tl]{\footnotesize Proceedings of the Ninth Prague Topological Symposium}}
\put(143,50){\makebox(0,0)[tl]{\footnotesize Contributed papers from the symposium held in}}
\put(143,34){\makebox(0,0)[tl]{\footnotesize Prague, Czech Republic, August 19--25, 2001}}
\end{picture}
\end{center}
\vspace{0.25in}
\setcounter{page}{155}
\title{Special metrics}
\author{Yasunao Hattori}
\address{Department of Mathematics, 
Shimane University, Matsue, Shimane, 690-8504 Japan}
\email{hattori@math.shimane-u.ac.jp}
\thanks{The author was an invited speaker at the Ninth Prague Topological Symposium.}
\thanks{Yasunao Hattori,
{\em Special metrics},
Proceedings of the Ninth Prague Topological Symposium, (Prague, 2001),
pp.~155--164, Topology Atlas, Toronto, 2002}
\begin{abstract}
This is a survey on special metrics.
We shall present some results and open questions on special metrics
mainly appeared in the last 10 years.
\end{abstract}
\maketitle

\section{Introduction}
This is a survey on special metrics.
We shall present some results and open questions on special metrics
mainly appeared in the last 10 years.
We shall consider {\it metrizable spaces} and all metrics induce the
original topology of a given metrizable space. 
For a metrizable space $X$, there are many metrics which induce the
original topology of $X$. 
Some of them may determine topological properties of $X$. 
For example, it is well known that a metrizable space $X$ is {\it
separable} if and only if $X$ admits a {\it totally bounded} metric, and
$X$ is {\it compact} if and only if $X$ admits a {\it complete totally
bounded} metric. 

\section{Properties related to midsets}

Let $(X, \rho)$ be a metric space and $y \ne z$ be distinct points of $X$. 
A set of the form 
$$M(y,z) = \{x \in X : \rho (x, y) = \rho(x,z)\}$$
is called a {\bf midset} or a {\bf bisector}.

Midsets are a geometrically intuitive concept, and several topological
properties can be approached through midsets. 

For example, the covering dimension of a separable metrizable space $X$
can be characterized by midsets: 

\begin{thm}[Janos-Martin \cite{JM}] 
A separable metrizable space $X$ has $\dim X \leq n$ if and only if $X$
admits a totally bounded metric $\rho$ such that $\dim M \leq n - 1$ for
every midset $M$ in $X$.
\end{thm}

Furthermore, metrics which have special midsets may determine the 
topological structures of spaces. 

\begin{defn} 
For a natural number $n$, a metric space $(X,\rho)$ is said to have the
{\bf $n$-points midset property}, abbreviated as the {\bf $n$-MP},
if every midset in $X$ has exactly $n$ points. 
\end{defn}

The 1-MP, 2-MP and 3-MP are sometimes called the {\bf unique midset
property}, {\bf double midset property} and {\bf triple midset property},
abbreviated as the {\bf UMP}, {\bf DMP} and {\bf TMP}, respectively.

The real line with the usual metric is an example of the space having the
UMP, and the circle in the two-dimensional Euclidean plane is an example
of the space having the DMP. 

Berard \cite{B} proved that a connected metric space with the UMP is
homeomorphic to an interval.
Further, Nadler \cite{Nadler} proved that every non-degenerate component 
of a metric space with the UMP is homeomorphic to an interval and that a
separable, locally compact metric space with the UMP is homeomorphic to a
subspace of the real line. 
However, the following question remains open.

\begin{ques}
Is a separable metric space with the UMP homeomorphic to a subspace of
the real line?
\end{ques} 

Concerning the question above, we have the following example.

\begin{exmp}[Hattori-Ohta \cite{HO}]
There is a separable metric space $(X, \rho)$ such that the cardinality
of any midset is at most 1, but $X$ is not suborderable, and hence $X$ is
not homeomorphic to a subspace of the real line.
\end{exmp}

Metrizable spaces that are homeomorphic to subspaces of the real line are
also
characterized as follows. 

\begin{thm}[Hattori-Ohta \cite{HO}] 
A separable metrizable space $X$ is homeomorphic to a subspace of the
real line if and only if $X$ admits a metric satisfying both of the
following conditions:
\begin{itemize} 
\item[(i)]
The cardinality of any midset is at most 1.
\item[(ii)]
The cardinality of any subset consisting of points which are equidistant
from a point is at most 2.
\end{itemize}
Furthermore, if $X$ is locally compact, then $X$ is homeomorphic to a 
subspace of the real line if and only if $X$ admits a metric satisfying
the condition (i) only.
\end{thm}

Now, we have an open question.

\begin{ques}[\cite{HO}]
Is every rim-compact (= each point has a neighborhood base consisting of
the open sets with compact boundaries), separable metric space 
$(X, \rho)$ satisfying the condition (i) homeomorphic to a subspace of the
real line?
\end{ques}

A metrizable space that admits a metric with $n$-points midset property
is also said to have the $n$-points midset property.
Several spaces are known to have the UMP, or do not have the UMP.

\begin{thm}[Ohta-Ono \cite{OO}, Ito-Ohta-Ono \cite{IOO}] 
\mbox{}
\begin{itemize}
\item[(a)]
Let $I$ and $J$ be separated intervals in the real line ${\mathbb R}$. 
Then $I \cup J$ admits a metric which has the UMP if and only if at least
one of $I$ and $J$ is not
compact.
\item[(b)]
The union of odd numbers of disjoint closed intervals in ${\mathbb R}$ admits a
metric which has the UMP.
\item[(c)]
The subsets $[0, 1] \cup {\mathbb Z}$ and $[0, 1] \cup {\mathbb Q}$ of ${\mathbb R}$ do
not admit metrics which have the UMP, where ${\mathbb Z}$ and ${\mathbb Q}$ denote
the sets of the integers and rational numbers, respectively.
\item[(d)]
Let $X$ be the union of at most countably many subsets $X_n$ of ${\mathbb R}$.
If each $X_n$ is either an interval or totally disconnected and if at
least one of $X_n$ is a noncompact interval, then $X$ admits a metric
which has the UMP.
\item[(e)]
A discrete space $D$ admits a metric which has the UMP if and only if 
$|D| \ne 2, 4$ and $|D| \leq {\mathfrak c}$, where $|D|$ denotes the
cardinality of $D$ and ${\mathfrak c}$ is the cardinality of the
continuum.
\item[(f)]
Let $D$ be the discrete space with $|D| \leq {\mathfrak c}$. 
Then the product of countably many copies of $D$ admits a metric which has
the UMP. In particular, the Cantor set and the space of irrationals admit
metrics which have the UMP.
\end{itemize}
\end{thm}

\begin{ques}[\cite{OO}]
Does every subspace of the real line ${\mathbb R}$ containing a noncompact
intervals as a clopen set have the UMP?
\end{ques}

\begin{ques}[\cite{IOO}]
Is there a subspace $X$ of ${\mathbb R}^n$ or the Hilbert space such that $X$
is homeomorphic to the Cantor set or the space of irrationals and the
metric inherited from the usual metric has the UMP?
\end{ques}

The UMP of a finite discrete space can be considered in terms of graph 
theory.
Let $G$ be a {\bf simple graph} (i.e., a graph which does not contain
neither multiple edges nor loops).
By a {\bf colouring} of $G$ we mean a map defined on the set of edges
$E(G)$ of $G$.
A coloring $\varphi$ of $G$ is said to have the {\bf unique midset 
property} if for every pair of distinct vertices $x$ and $y$ there is a
unique vertex $p$ such that $xp$ and $yp$ are edges of $G$ and
$\varphi(xp) = \varphi(yp)$.

\begin{thm}[Ito-Ohta-Ono \cite{IOO}] 
Let $K_n$ be the complete graph (i.e., each vertex of $K_n$ is adjacent to
every other vertices) with $n$ vertices.
Then, a finite discrete space with $n$ points has the UMP if and only if 
there is a colouring $\varphi$ of the complete graph $K_n$ with the UMP.
\end{thm}

Let $G$ be a finite graph having a colouring with the UMP.
Then we denote by ump$(G)$ the smallest number of colours required for a
colouring of $G$ with the UMP; i.e.,
$$
\ump(G) = 
\min \{|\varphi(E(G))| : 
\mbox{ $\varphi$ is a colouring of $G$ with the UMP}\}.
$$

\begin{thm}[Ito-Ohta-Ono \cite{IOO}] 
For a complete graph $K_n$, we have the following.
\begin{itemize}
\item[(a)]
For each $n \geq 0$, ump$(K_{2n+1}) = n$.
\item[(b)]
For each $n \geq 3$, ump$(K_{2n}) \leq 2n-1$. 
\end{itemize}
In particular, we have 
\begin{itemize}
\item $\ump(K_6) = 4$, 
\item $\ump(K_8) = 5$,
\item $\ump(K_{10}) = 5$,
\item $\ump(K_{12}) \leq 8$, 
\item $\ump(K_{14}) \leq 10$.
\end{itemize}
\end{thm}

\begin{ques}[\cite{IOO}]
Determine the values of ump$(K_{2n})$ for each $n \geq 6$.
\end{ques} 

The following {\it double midset conjecture} seems to be very interesting.

\begin{dmc}
A continuum (= non-degenerate connected compact metric space) having the
DMP must be a simple closed curve. 
\end{dmc}

The conjecture still remains open. 
However, several partial results about the conjecture are known.

\begin{thm}
A continuum $X$ with the DMP satisfying the either of the following
conditions is a simple closed curve:
\begin{itemize}
\item[(a)]
{\rm (L. D. Loveland and S. G. Wayment \cite{LW})}
$X$ contains a continuum with no cut points.
\item[(b)]
{\rm (L. D. Loveland \cite{L76})}
The midset function $M : \{(x,y) \in X \times X : x \ne y\} \to 2^X$ is
continuous.
\item[(c)]
{\rm (L. D. Loveland \cite{L91})}
$X \subset {\mathbb R}^2$ with the Euclidean metric.
\end{itemize}
\end{thm}

Furthermore, L. D. Loveland and S. M. Loveland proved: 

\begin{thm}[L. D. Loveland and S. M. Loveland \cite{LL}] 
Every continuum in the Euclidean plane with the $n$-MP for $n \geq 1$ must either be a simple closed curve or an arc.
\end{thm}

We shall consider a more general setting of the double midset conjecture.

\begin{smc}
A non-degenerate compact metric space such that all of its midsets are
homeomorphic to an $(n-1)$-sphere $S^{n-1}$ is homeomorphic to an
$n$-sphere $S^{n}$. 
\end{smc}

The double midset conjecture can be considered as the $0$-sphere midset
conjecture. 
A space is said to have the {\bf $k$-sphere midset property} if each
midset of $X$ is homeomorphic to a $k$-sphere $S^k$. 
We have a few results in this direction: 

\begin{thm}
\mbox{}
\begin{itemize}
\item[(a)] {\rm (L. D. Loveland \cite{L77})}
If $X$ is a metric space with 1-sphere midset property, and if $X$
contains a subset homeomorphic to 2-sphere, then $X$ is a 2-sphere. 
\item[(b)] {\rm (L. D. Loveland \cite{L77})}
If $X$ is a nondegenerate compact metric space with 1-sphere midset
property and every simple closed curve separates $X$, then $X$ is a 2-sphere. 
\item[(c)] {\rm (W. D\c{e}bski, K. Kawamura and K. Yamada \cite{DKY})}
Let $n \geq 3$ and $X$ be a nondegenerate compact subset of the
$n$-dimensional Euclidean space such that each of its midsets is a convex
$(n-2)$-sphere (= the boundary of a convex $(n-1)$-cell). Then $X$ is a
convex $(n-1)$-sphere.
\end{itemize}
\end{thm}

\section{Metrics having special properties}

There are several topology preserving metrics in a metrizable space $X$. 
Some of them may have special properties.
From this point of view, we have the following results:
(But, the results below are also motivated by some results concerning
special metrics characterizing topological dimension.) 

\begin{thm}
Every metrizable space $X$ admits the following metrics $\rho_1$,
$\rho_2$, $\rho_3$ and $\rho_4$.
\begin{itemize}
\item[(a)] {\rm (Nagata \cite{Nagata1} or see \cite{HN})}
For every $x \in X$ and each sequence $\{y_1, y_{2}, \dots\}$ in $X$
there is a triplet of indices $i$, $j$ and $k$ such that $i \ne j$ and
$\rho_1(y_i, y_j) \leq \rho_1(x, y_k)$.
\item[(b)] {\rm (Nagata, see \cite{NagataB})}
$\{S_{\varepsilon}(x) : x \in X\}$ is closure-preserving for each
$\varepsilon > 0$, where 
$S_{\varepsilon}(x) = \{y \in X : \rho_2(x,y) < \varepsilon\}$ is an
$\varepsilon$-ball with respect to $\rho_2$.
\item[(c)] {\rm (Hattori \cite{H})} 
For every $x \in X$ and every sequence $\{y_1, y_{2}, \dots\}$ in $X$ 
with $\rho_3(x, y_i) \geq \delta$ for infinitely many $y_i$ for some
$\delta > 0$, there is a pair of distinct indices $i$ and $j$ such that
$\rho_3(y_i, y_j) \leq \rho_3(x, y_i)$. 
\item[(d)] {\rm (Hattori \cite{H})}
$X$ has a $\sigma$-locally finite open base consisting of open balls with
respect to $\rho_4$.
\end{itemize}
\end{thm}

Concerning (c) and (d) in the theorem above, the following are still
remain open.

\begin{ques}[\cite{Nagata1} or see \cite{HN}]
Can we drop the condition ``$\rho(x, y_i) \geq \delta$ for infinitely many
$y_i$ for some $\delta > 0$'' in (c)?
I.e., does every metrizable space admit a metric $\rho$ such that for
every $x \in X$ and every sequence $\{y_1, y_{2}, \dots \}$ in $X$, there
is a pair of distinct indices $i$ and $j$ such that 
$\rho_(y_i, y_j) \leq \rho_(x, y_i)$?
\end{ques} 

\begin{ques}[\cite{Nagata1}]
Does every metrizable space $X$ admit a metric $\rho$ such that 
$\mathcal B_n = \{S_{1/n}(x) : x \in X_n\}$ is discrete in $X$ for some 
$X_n \subset X$, $n = 1, 2, \dots$ and 
$\mathcal B = \bigcup_{n=1}^{\infty}\mathcal B_n$ is a base for $X$?
\end{ques} 

Nagata \cite{Nagata2} also asked whether if every metrizable space admits
a metric such that $\{S_{\varepsilon}(x) : x \in X\}$ is hereditarily
closure-preserving for each $\varepsilon > 0$. This is answered
negatively by Ziqiu-Junnila and Balogh-Gruenhage independently as follows.

\begin{thm}[Ziqiu-Junnila \cite{ZJ}]
The hedgehog space of the weight $\geq 2^{{\mathfrak c}^+}$ does
not admit a metric such that $\{S_{\varepsilon}(x) : x \in X\}$ is
hereditarily closure-preserving for each $\varepsilon > 0$. 
\end{thm}

\begin{thm}[Balogh-Gruenhage \cite{BG}] 
A metrizable space $X$ admits a metric such that
$\{S_{\varepsilon}(x) : x \in X\}$ is hereditarily closure-preserving for
each $\varepsilon > 0$ if and only if $X$ is strongly metrizable, where a
metrizable space $X$ is called strongly metrizable if there is a base for
$X$ which is the union of countably many star-finite open coverings.
In particular, the hedgehog space of weight $\omega_1$ does not have such
metric.
\end{thm}

We notice that there is another chracterization of strongly metrizable
spaces by a special metric (Hattori \cite{H}) : A metrizable space $X$ is
strongly metrizable if and only if $X$ admits a metric $\rho$ such that
for every $\varepsilon >0$, every $x \in X$ and every sequence $\{y_1,
y_2, \dots \}$ in $X$ with $\rho(S_{\varepsilon}(x) , y_i) < \varepsilon$
for all $i$, there is a pair of distinct indicies $i$ and $j$ such that
$\rho(y_i, y_j) < \varepsilon$.

The following questions asked by Nagata in \cite{Nagata1} seems to be open.

\begin{ques}
If $X$ is a rim-separable (= every point has a neighborhood base
consisting of open sets with separable boundaries) metrizable space, then
does $X$ admit a metric $\rho$ such that $\Bd S_{\varepsilon}(x)$ is
separable for each $\varepsilon > 0$ and $x \in X$?
\end{ques}

\begin{ques}
If $X$ is a rim-compact (= every point has a neighborhood base consisting
of open sets with compact boundaries) metrizable space, then does $X$
admit a metric $\rho$ such that $\Bd S_{\varepsilon}(x)$ is compact for
each $\varepsilon > 0$ and $x \in X$?
\end{ques}

\section{Ultrametric spaces}

A metric $\rho$ is called an {\bf ultrametric} (or a {\bf non-Archimedean
metric}) if $\rho$ satisfies the strong triangle inequality: 
$$\max\{\rho (x, y), \rho (y,z)\} \geq 
\rho (x,z) \mbox{ for all $x, y, z \in X$}.$$

Ultrametric spaces are topologically characterized independently by F.
Hausdorff, J. de Groot (\cite{Groot} and also see \cite[Problem 4.1.H]{E})
and K. Morita \cite{M}: A metrizable space $X$ admits an ultrametric if
and only if $\dim X = 0$, where $\dim X$ denotes the covering dimension of
$X$. The characterization was extended to higher dimension by J. Nagata
\cite{Nagata} and P. Ostrand \cite{Ostrand}. Ultrametric spaces are also
studied in dimension theory from more geometrical points of view.

Let ${\mathbb R}_+$ denote the set of non-negative real numbers and $C$, 
$s \in {\mathbb R}_+$. 

\begin{defn}
A metric space $(X, \rho)$ is {\bf $(C, s)$-} {\bf homogeneous} if the
inequality $|X_0| \leq C(b/a)^s$ holds for $a > 0$, $b > 0$ and
$X_0 \subset X$ provided that $b \geq a$ and that 
$a \leq \rho(x, y) \leq b$ holds for every pair of distinct points $x$ and
$y$ of $X_0$. 

The space $(X, \rho)$ is said to be {\bf $s$-homogeneous} if it is 
$(C, s)$-homogeneous for some $C \in {\mathbb R}_+$. 
\end{defn}

\begin{defn}
We define the {\bf Assouad dimension} $\dim_A (X, \rho)$ of a metric
space $(X, \rho)$ as follows:
$$\dim_A (X, \rho) = 
\inf\{s \in {\mathbb R}_+ : (X, \rho) \mbox{ is $s$-homogeneous}\},$$ 
if the infimum exists. 
Otherwise, we define $\dim_A (X, \rho) = \infty$. 
\end{defn}

Then we have the following theorem.

\begin{thm}
An ultrametric space $(X, \rho)$ can be bi-Lipschitz embedded in the
$n$-dimensional Euclidean space ${\mathbb R}^n$ if and only if 
$\dim_A (X, \rho) < n$, where a mapping $f : X \to {\mathbb R}^n$ is said to be
a {\bf bi-Lipschitz embedding} if there exists a real number 
$\alpha \geq 1$ such that for every $x$, $y \in X$, 
$\frac1{\alpha}\rho(x, y) \leq \|f(x) - f(y)\| \leq \alpha\rho(x, y)$
holds. 
\end{thm}

The ``if'' part of the theorem is proved by Luukkainen and
Movahedi-Lankarani \cite{LML} and the ``only if'' part is proved by Luosto
\cite{Luosto}.

Several metrically universal properties of ultrametric spaces are known. 
A. J. Lemin (\cite{Lemin1}) proved that every ultrametric space of weight 
$\tau$ is isometrically embedded in the generalized Hilbert space
$H^{\tau}$. 
Furthermore, A. J. Lemin and V. A. Lemin proved the following:

\begin{thm}[A. J. Lemin and V. A. Lemin \cite{Lemin2}] 
For every cardinal $\tau$ there is an ultrametric space 
$(LW_{\tau}, \rho)$ such that every ultrametric space of weight 
$ \leq \tau$ is isometrically embedded in $(LW_{\tau}, \rho)$. 
\end{thm}

We notice that the weight of the space $(LW_{\tau}, \rho)$ is
$\tau^{\aleph_0}$ and this is the best possible if we consider the
countable case. In fact, if an ultrametric space $(X, \rho)$ isometrically
contains all two-point ultrametric spaces, then the weight of 
$X \geq {\mathfrak c}$, where ${\mathfrak c}$ denotes the cardinality of
the continuum. 
On the other hand, J. Vaughan \cite{V} proved the following.

\begin{thm}[J.~Vaughan \cite{V}] 
Under the assumption of the singular cardinal hypothesis, for every
cardinal $\tau$ satisfying ${\mathfrak c} < \tau < \tau^{\omega}$
there exists an ultrametric space $(LW_{\tau}', \rho')$ such that the
weight of $(LW_{\tau}', \rho')$ is $\tau$ and every ultrametric space of
weight $\leq \tau$ is isometrically embedded in $(LW_{\tau}', \rho')$. 
\end{thm}

In the case of spaces consisting of finite points, we have the following.

\begin{thm}[A. J. Lemin and V. A. Lemin \cite{Lemin2}]
Every ultrametric space consisting of $n+1$ points is isometrically
embedded in the $n$-dimensional Euclidean space ${\mathbb R}^n$ and there is no
ultrametric space $X$ consisting of $n+1$ points such that $X$ is
isometrically embedded in the $k$-dimensional Euclidean space ${\mathbb R}^k$
for $k < n$.
\end{thm}

Recently, the theory of ultrametric spaces has been developped in several
branches of mathematics as well as physics, biology and information
sciences. In particular, the theory of ultrametric spaces is applied to
domain theory and logic programming.

We shall introduce some notion and terminology from the domain theory. 
Let $(P, \leq)$ be a partially ordered set.
A partial order $\leq$ on $(P, \leq)$ is said to be {\bf directed 
complete} if every directed subset $D$ of $P$ has a least upper bound
$\bigvee D$ and a partially ordered set $(P, \leq)$ is said to be a {\bf
directed complete partially ordered set} if $\leq$ is directed complete.
(A directed complete partially ordered set is often abbreviated as a 
{\bf dcpo}.)
Let $P$ be a dcpo and $x, y \in P$.
We say $x$ is {\bf way below} $y$, written by $x \ll y$, if for every
directed subset $D$ of $P$ with $y \leq \bigvee D$ there is $d \in D$
such that $x \leq d$.
An element $x$ of a dcpo $P$ is {\bf compact} if $x \ll x$ and we denote
the set of compact elements of $P$ by $P_c$.
A dcpo $P$ is called {\bf algebraic} ({\bf $\omega$-algebraic}) if ($P_c$
is countably infinite and) every element of $P$ is a directed sup of
compact elements, i.e., for each element $x \in P$, 
$\{y \in P_c : y \leq x\}$ is directed and 
$x = \bigvee \{y \in P_c : y \leq x\}$.
An algebraic dcpo is called a {\bf domain}.
We call that a dcpo $P$ is {\bf $\omega$-continuous} if there is a
countably infinite subset $B$ of $P$ such that for every $x \in P$, 
$\{y : y \in B, \mbox{ $y \ll x$} \}$ is directed and 
$\bigvee \{y : y \in B, \mbox{ $y \ll x$} \} = x$.

\begin{defn}
The {\bf Scott topology} of a partially ordered set $P$ is defined as
follows: A subset $U \subset P$ is open if 
\begin{itemize}
\item[(i)]
for every $x \in U$, $\{ y \in P : y \geq x\} \subset U$, and 
\item[(ii)]
for every directed subset $D$ in $P$ with $\bigvee D \in U$, 
$D \cap U \ne \emptyset$. 
\end{itemize}
\end{defn}

If $P$ is a domain, then the Scott topology is generated by the subbase
consisting of sets of the form $\{y \in P : y \geq x\}$ for $x \in P_c$.

\begin{defn}
The {\bf Lawson topology} of a partially ordered set $P$ is the supremun
of the Scott topology and the weak$^d$ topology, where the {\bf weak$^d$
topology} is the topology determined by the closed base consisting of
sets of the form $\{y \in P : y \geq x\}$ for $x \in P$.
\end{defn}

\begin{defn}
Let $\Max (P)$ be the set of maximal elements of $P$. 
A {\bf computational model} for a topological space $X$ is an
$\omega$-continuous diretected complete partially ordered set $P$ with an 
embedding $i : X \to \Max (P)$ which satisfy the following
conditions:
\begin{itemize}
\item[(1)] 
The restrictions of the Scott topology and the Lawson topology to $\Max
(P)$ coincide. 
\item[(2)]
The embedding $i : X \to \Max (P)$ is a homeomorphism. 
\end{itemize}
\end{defn}

Then, B. Flagg and R. Kopperman (\cite{FK}) proved the following.

\begin{thm} 
A topological space $X$ has an $\omega$-algebraic computational model if
and only if $X$ is a complete separable ultrametric space. 
\end{thm}

We refer the reader to \cite{LM} for terminology and recent developments
of domain theory related to topology. 
We also refer the reader to \cite{Lemin2} for a brief historical 
introduction to ultrametric spaces.

\subsection*{Acknowledgements}
The author would like to his thanks to the referee for his (her) helpful 
comments.

%\bibliographystyle{amsplain}
%\bibliography{15}

\begin{thebibliography}{10}

\bibitem{BG}
Zoltan Balogh and Gary Gruenhage, \emph{When the collection of $\epsilon$-balls
  is locally finite}, To appear in Topology Appl.

\bibitem{B}
Anthony~D. Berard, Jr., \emph{Characterizations of metric spaces by the use of
  their midsets: {I}ntervals}, Fund. Math. \textbf{73} (1971/72), no.~1, 1--7.
  \MR{45 \#4368}

\bibitem{Groot}
J.~de~Groot, \emph{Non-archimedean metrics in topology}, Proc. Amer. Math. Soc.
  \textbf{7} (1956), 948--953. \MR{18,325a}

\bibitem{DKY}
W.~D{\c{e}}bski, K.~Kawamura, and K.~Yamada, \emph{Subsets of $\mathbb{R}^n$
  with convex midsets}, Topology Appl. \textbf{60} (1994), no.~2, 109--115.
  \MR{95j:52003}

\bibitem{E}
Ryszard Engelking, \emph{Theory of dimensions finite and infinite}, Heldermann
  Verlag, Lemgo, 1995. \MR{97j:54033}

\bibitem{FK}
Bob Flagg and Ralph Kopperman, \emph{Computational models for ultrametric
  spaces}, Mathematical foundations of programming semantics (Pittsburgh, PA,
  1997), Elsevier, Amsterdam, 1997, p.~9 pp. (electronic). \MR{98h:68151}

\bibitem{H}
Yasunao Hattori, \emph{On special metrics characterizing topological
  properties}, Fund. Math. \textbf{126} (1986), no.~2, 133--145. \MR{87m:54103}

\bibitem{HN}
Yasunao Hattori and Jun-iti Nagata, \emph{Special metrics}, Recent progress in
  general topology (Prague, 1991), North-Holland, Amsterdam, 1992,
  pp.~353--367. \MR{1 229 131}

\bibitem{HO}
Yasunao Hattori and Haruto Ohta, \emph{A metric characterization of a subspace
  of the real line}, Topology Proc. \textbf{18} (1993), 75--87. \MR{95k:54052}

\bibitem{IOO}
Munehiko It{\=o}, Haruto Ohta, and Jin Ono, \emph{A graph-theoretic approach to
  the unique midset property of metric spaces}, J. London Math. Soc. (2)
  \textbf{60} (1999), no.~2, 353--365. \MR{2000i:05176}

\bibitem{JM}
Ludvik Janos and Harold Martin, \emph{Metric characterizations of dimension for
  separable metric spaces}, Proc. Amer. Math. Soc. \textbf{70} (1978), no.~2,
  209--212. \MR{57 \#13876}

\bibitem{LM}
Jimmie~D. Lawson and Michael Mislove, \emph{Problems in domain theory and
  topology}, Open problems in topology, North-Holland, Amsterdam, 1990,
  pp.~349--372. \MR{1 078 658}

\bibitem{Lemin1}
A.~Yu. Lemin, \emph{Isometric imbedding of isosceles (non-{A}rchimedean) spaces
  into {E}uclidean ones}, Dokl. Akad. Nauk SSSR \textbf{285} (1985), no.~3,
  558--562, English translation: Soviet Math. Dokl. 32 (1985), no. 3, 740--744.
  \MR{87h:54056}

\bibitem{Lemin2}
Alex~J. Lemin and Vladimir~A. Lemin, \emph{On a universal ultrametric space},
  Topology Appl. \textbf{103} (2000), no.~3, 339--345. \MR{2001b:54035}

\bibitem{L76}
L.~D. Loveland, \emph{A metric characterization of a simple closed curve},
  General Topology and Appl. \textbf{6} (1976), no.~3, 309--313. \MR{53 \#9170}

\bibitem{L77}
\bysame, \emph{Metric spaces with connected midsets}, Houston J. Math.
  \textbf{3} (1977), no.~4, 495--501. \MR{57 \#17601}

\bibitem{L91}
\bysame, \emph{The double midset conjecture for continua in the plane},
  Topology Appl. \textbf{40} (1991), no.~2, 117--129. \MR{92h:54047}

\bibitem{LL}
L.~D. Loveland and S.~M. Loveland, \emph{Equidistant sets in plane triodic
  continua}, Proc. Amer. Math. Soc. \textbf{115} (1992), no.~2, 553--562.
  \MR{92i:54034}

\bibitem{LW}
L.~D. Loveland and S.~G. Wayment, \emph{Characterizing a curve with the double
  midset property}, Amer. Math. Monthly \textbf{81} (1974), 1003--1006. \MR{54
  \#6103}

\bibitem{Luosto}
Kerkko Luosto, \emph{Ultrametric spaces bi-{L}ipschitz embeddable in
  $\mathbb{R}^n$}, Fund. Math. \textbf{150} (1996), no.~1, 25--42.
  \MR{97j:54030}

\bibitem{LML}
Jouni Luukkainen and Hossein Movahedi-Lankarani, \emph{Minimal bi-{L}ipschitz
  embedding dimension of ultrametric spaces}, Fund. Math. \textbf{144} (1994),
  no.~2, 181--193. \MR{95i:54031}

\bibitem{M}
Kiiti Morita, \emph{Normal families and dimension theory for metric spaces},
  Math. Ann. \textbf{128} (1954), 350--362. \MR{16,501h}

\bibitem{Nadler}
Sam~B. Nadler, Jr., \emph{An embedding theorem for certain spaces with an
  equidistant property}, Proc. Amer. Math. Soc. \textbf{59} (1976), no.~1,
  179--183. \MR{53 \#14433}

\bibitem{Nagata}
Jun-iti Nagata, \emph{On a special metric and dimension}, Fund. Math.
  \textbf{55} (1964), 181--194. \MR{31 \#5186}

\bibitem{Nagata1}
\bysame, \emph{A survey of special metrics}, Proc. General Topology Symposium,
  Okayama Univ., 1982, pp.~36--45.

\bibitem{NagataB}
\bysame, \emph{Modern dimension theory}, Heldermann Verlag, Berlin, 1983.
  \MR{84h:54033}

\bibitem{Nagata2}
\bysame, \emph{Remarks on metrizability and generalized metric spaces},
  Topology Appl. \textbf{91} (1999), no.~1, 71--77. \MR{2000c:54022}

\bibitem{OO}
Haruto Ohta and Jin Ono, \emph{The unique midpoint property of a subspace of
  the real line}, Topology Appl. \textbf{104} (2000), no.~1-3, 215--226.
  \MR{2001e:54058}

\bibitem{Ostrand}
Phillip~A. Ostrand, \emph{A conjecture of {J}. {N}agata on dimension and
  metrization}, Bull. Amer. Math. Soc. \textbf{71} (1965), 623--625. \MR{31
  \#1655}

\bibitem{V}
Jerry~E. Vaughan, \emph{Universal ultrametric spaces of smallest weight},
  Topology Proc. \textbf{24} (1999), no.~Summer, 611--619 (2001). \MR{1 876
  391}

\bibitem{ZJ}
Ziqiu Y. and H.~Junnila, \emph{On a special metric}, To appear in Houston J.
  Math.

\end{thebibliography}
\providecommand{\bysame}{\leavevmode\hbox to3em{\hrulefill}\thinspace}
\providecommand{\MR}{\relax\ifhmode\unskip\space\fi MR }
\providecommand{\MRhref}[2]{%
  \href{http://www.ams.org/mathscinet-getitem?mr=#1}{#2}
}
\providecommand{\href}[2]{#2}

\end{document}